\documentclass[a4paper,11pt]{amsart}

\usepackage{graphicx}
\usepackage{mathptmx}
\usepackage{amsmath}
\usepackage{amssymb}
\usepackage{enumitem}
\usepackage{xcolor}

\newmuskip\pFqmuskip

\newcommand*\pFq[6][8]{%
  \begingroup 
  \pFqmuskip=#1mu\relax
  \mathcode`=\string"8000
  \begingroup\lccode`\~=`\,
  \lowercase{\endgroup\let~}\pFqcomma
  F^{#2}_{#3}{\left(\genfrac..{0pt}{}{#4}{#5}\bigg|#6\right)}%
  \endgroup
}
\newcommand{\pFqcomma}{\mskip\pFqmuskip}

\newtheorem{theorem}{Theorem}[section]

\begin{document}

\title[Spivey's type recurrence relation for degenerate Bell polynomials]{Spivey's type recurrence relation for degenerate Bell polynomials}

\author{Taekyun  Kim}
\address{Department of Mathematics, Kwangwoon University, Seoul 139-701, Republic of Korea}
\email{tkkim@kw.ac.kr}
\author{Dae San  Kim}
\address{Department of Mathematics, Sogang University, Seoul 121-742, Republic of Korea}
\email{dskim@sogang.ac.kr}


\begin{abstract}
The aim of this paper is to derive a recurrence relation for the degenerate Bell polynomials by using the operators $X$ and $D$ satisfying the commutation relation $DX-XD=1$. Here $X$ is the `multiplication by $x$' operator and $D=\frac{d}{dx}$. This recurrence relation is a generalization of Spivey's recurrence relation for the Bell numbers. We also obtain a recurrence relation for the degenerate $r$-Bell polynomials by using the same operators.
\end{abstract}

\maketitle

\section{Introduction}
The Stirling number of the second kind ${n \brace k}$ enumerates the number of partitiions of the set $[n]=\{1,2,\dots,n\}$ into $k$ nonempty disjoint subsets (see [15,27]). The sum of the Stirling numbers of the second kind ${n \brace k}$, called the Bell number and denoted by $\phi_{n}$, counts the total number of partitions of a set with $n$ elements (see [11,12]). That is, $\phi_{n}=\sum_{k=0}^{n}{n \brace k}$. The Bell polynomial $\phi_{n}(x)=\sum_{k=0}^{n}{n \brace k}x^{k}$ is a natural polynomial extension of the Bell number $\phi_{n}$ (see [1,31]). \par
Let $r$ be a nonnegative integer. The $r$-Stirling number of the second kind ${n+r \brace k+r}_{r}$ enumerates the number of partitiions of $[n]$ into $k$ nonempty disjoint subsets in such a way that $1,2,\dots,r$ are in distinct nonempty disjoint subsets (see [6]). The sum of the $r$-Stirling numbers of the second kind ${n+r \brace k+r}_{r}$, called the $r$-Bell number and denoted by $\phi_{n}^{(r)}$, counts the total number of partitions of $[n]$ in such a way that $1,2,\dots,r$ are in distinct nonempty disjoint subsets. That is, $\phi_{n}^{(r)}=\sum_{k=0}^{n}{n+r \brace k+r}_{r}$. The $r$-Bell polynomial $\phi_{n}^{(r)}(x)=\sum_{k=0}^{n}{n+r \brace k+r}_{r}x^{k}$ is a natural polynomial extension of the $r$-Bell number $\phi_{n}^{(r)}$. \par
In recent years, many degenerate versions have been explored for special polynomials and numbers (see [7,9,13,17,19,20,21,23,25,28,29]), gamma functions (see [24]) and even umbral calculus (see [16]) by employing diverse tools, including combinatorial methods, generating functions, $p$-adic analysis, umbral calculus, probability theory, quantum mechanics, and operator theory. \par
In this paper, we consider degenerate versions of the Bell polynomials and the $r$-Bell polynomials, namely the degenerate Bell polynomials $\phi_{n,\lambda}(x)$ and the degenerate $r$-Bell polynomials $\phi_{n,\lambda}^{(r)}(x)$. \par
Spivey obtained the following recurrence relation for the Bell numbers:
\begin{equation}
\phi_{m+n}=\sum_{j=0}^{m}\sum_{k=0}^{n}\binom{n}{k}{m \brace j}j^{n-k}\phi_{k},\quad (m,n\ge 0),\quad (\mathrm{see}\ [22,35]).\label{0}	
\end{equation}
The aim of this paper is to derive Spivey's type recurrence relations for $\phi_{n,\lambda}(x)$ and $\phi_{n,\lambda}^{(r)}(x)$ by using the `multiplication by $x$' operator $X$ and the differentiation operator $D=\frac{d}{dx}$ (see [4,5,14,17,19,25,26,33]), satisfying the commutation relation $DX-XD=1$. Indeed, for the degenerate $r$-Bell polynomials $\phi_{n,\lambda}^{(r)}(x)$, we show the following:
\begin{equation*}
\phi_{n+m,\lambda}^{(r)}(x)=\sum_{k=0}^{m}\sum_{l=0}^{n}\binom{n}{l}{m+r \brace k+r}_{r,\lambda}(k-m\lambda)_{n-l,\lambda}x^{k}\phi_{l,\lambda}^{(r)}(x).
\end{equation*} \par
As general references, we let the reader refer to [2,3,8,10,34,36]. In the rest of this section, we recall the necessary facts that are needed throughout this paper.

\vspace{0.1in}

For any nonzero $\lambda\in\mathbb{R}$, the degenerate exponentials are defined by Kim-Kim as
\begin{equation}
e_{\lambda}^{x}(t)=\sum_{n=0}^{\infty}(x)_{n,\lambda}\frac{t^{n}}{n!},\quad\mathrm{and}\quad e_{\lambda}(t)=e_{\lambda}^{1}(t),\quad (\mathrm{see}\ [20]),
\end{equation}
where
\begin{equation}
(x)_{0,\lambda}=1,\qquad (x)_{n,\lambda}=x(x-\lambda)(x-2\lambda)\cdots\big(x-(n-1)\lambda\big),\quad (n\ge 1).\label{2}
\end{equation}
Note that
\begin{equation*}
\lim_{\lambda\rightarrow 0}e_{\lambda}^{x}(t)=e^{xt}.
\end{equation*} \par
The Stirling numbers of the second kind are defined by
\begin{equation}
x^{n}=\sum_{k=0}^{n}{n \brace k}(x)_{k},\quad (n\ge 0),\quad (\mathrm{see}\ [15,16,27]), \label{3}
\end{equation}
where
\begin{displaymath}
(x)_{0}=1,\quad (x)_{n}=x(x-1)(x-2)\cdots(x-n+1),\quad (n\ge 1).
\end{displaymath}
The degenerate Stirling numbers of the second kind are given by
\begin{equation}
(x)_{n,\lambda}=\sum_{k=0}^{n}{n \brace k}_{\lambda}(x)_{k},\quad (n\ge 0),\quad (\mathrm{see}\ [20,23,25,29,32]).\label{4}
\end{equation}
From \eqref{4}, we note that
\begin{equation}
\frac{1}{k!}\big(e_{\lambda}(t)-1\big)^{k}=\sum_{n=k}^{\infty}{n \brace k}_{\lambda}\frac{t^{n}}{n!}, \label{5}	
\end{equation}
where $k$ is a nonnegative integer. \par
For $r\in\mathbb{Z}$ with $r\ge 0$, the degenerate $r$-Stirling numbers of the second kind are defined by
\begin{equation}
(x+r)_{n,\lambda}=\sum_{k=0}^{n}{n+r \brace k+r}_{r,\lambda}(x)_{k},\quad (n\ge 0),\quad (\mathrm{see}\ [19]).\label{6}
\end{equation}
By \eqref{6}, we get
\begin{equation}
\frac{1}{k!}\big(e_{\lambda}(t)-1\big)^{k}e_{\lambda}^{r}(t)=\sum_{n=k}^{\infty}{n+r \brace k+r}_{r,\lambda}\frac{t^{n}}{n!},\quad (k\ge 0). \label{7}
\end{equation} \par
Also, the Bell polynomials are given by
\begin{equation}
e^{x(e^{t}-1)}=\sum_{n=0}^{\infty}\phi_{n}(x)\frac{t^{n}}{n!},\quad (\mathrm{see}\ [1,31]). \label{8}
\end{equation}
We note that $\phi_{n}(1)=\phi_{n},\ (n\ge 0)$, are the Bell numbers (see [11,12]) and
\begin{displaymath}
\phi_{n}(x)=\sum_{k=0}^{n}{n \brace k}x^{k}.
\end{displaymath}
In view of \eqref{8}, the degenerate Bell polynomials are defined by
\begin{equation}
e^{x(e_{\lambda}(t)-1)}=\sum_{n=0}^{\infty}\phi_{n,\lambda}(x)\frac{t^{n}}{n!},\quad (\mathrm{see}\ [18,30]). \label{9}
\end{equation}
By \eqref{5} and \eqref{9}, we get
\begin{equation}
\phi_{n,\lambda}(x)=\sum_{k=0}^{n}{n \brace k}_{\lambda}x^{k},\quad (n,k\ge 0).\label{10}
\end{equation}
Here $\phi_{n,\lambda}=\phi_{n,\lambda}(1),\ (n\ge 0)$, are called the degenerate Bell numbers. \par
The degenerate $r$-Bell polynomials are given by
\begin{equation}
e^{x(e_{\lambda}(t)-1)}e_{\lambda}^{r}(t)=\sum_{n=0}^{\infty}\phi_{n,\lambda}^{(r)}(x)\frac{t^{n}}{n!},\quad (\mathrm{see}\ [17]). \label{11}	
\end{equation}
By \eqref{7} and \eqref{11}, we get
\begin{equation}
\phi_{n,\lambda}^{(r)}(x)=\sum_{k=0}^{n}{n+r \brace k+r}_{r,\lambda}x^{k},\quad (n\ge 0).\label{12}	
\end{equation}

\section{Spivey's type recurrence relations for degenerate Bell and degenerate $r$-Bell polynomials}
Now, we define the operators $X$ and $D$ by
\begin{equation}
Xf(x)=xf(x)\quad\mathrm{and}\quad Df(x)=\frac{d}{dx}f(x)\label{14}	
\end{equation}
that satisfy $DX-XD=1$. \par
For $k\in\mathbb{N}$, we have
\begin{align}
DX^{k}-X^{k}D&=\big(DX^{k-1}-X^{k-1}D\big)X+X^{k-1}\big(DX-XD) \label{15}\\
 &=\big(DX^{k-1}-X^{k-1}D)X+X^{k-1}\nonumber\\
&=\big(DX^{k-2}-X^{k-2}D\big)X^{2}+2X^{k-1} \nonumber\\
&=\cdots \nonumber\\
&=\big(DX-XD\big)X^{k-1}+(k-1)X^{k-1}=kX^{k-1}.\nonumber
\end{align}
From \eqref{15}, we note that
\begin{equation}
\begin{aligned}
\big(XD\big)X^{k}&=X\big(DX^{k}\big)=X\big(X^{k}D+kX^{k-1}\big) \\
&=X^{k+1}D+kX^{k}=X^{k}\big(XD+k\big).
\end{aligned}	\label{16}
\end{equation}
By \eqref{16}, we get
\begin{equation}
\big(XD\big)X^{k}=X^{k}\big(XD+k\big),\quad (k\in\mathbb{N}).\label{17}	
\end{equation} \par
For $f(x)=x^{m}$, with $m \ge 0$, we have
\begin{equation}
\begin{aligned}
\big(XD\big)_{n,\lambda}f(x)&=\big(XD\big)_{n,\lambda}x^{m}=(m)_{n,\lambda}x^{m}=\sum_{k=0}^{n}{n \brace k}_{\lambda}(m)_{k}x^{m} \\
&=\sum_{k=0}^{n}{n \brace k}_{\lambda}X^{k}D^{k}x^{m}=\sum_{k=0}^{n}{n \brace k}_{\lambda}X^{k}D^{k}f(x).
\end{aligned}	\label{18}
\end{equation}
Thus, by comparing the coefficients on both sides of \eqref{18}, we get
\begin{equation}
\big(XD\big)_{n,\lambda}=\sum_{k=0}^{n}{n \brace k}_{\lambda}X^{k}D^{k},\quad (n\ge 0), 	\label{19}
\end{equation}
where
\begin{displaymath}
{n+1 \brace k}_{\lambda}={n \brace k-1}_{\lambda}+(k-n\lambda){n \brace k}_{\lambda},\quad (n\ge k\ge 1).
\end{displaymath}
From \eqref{10} and \eqref{19}, we note that
\begin{align}
\frac{1}{e^{x}}\big(XD\big)_{n,\lambda}e^{x}&=\frac{1}{e^{x}}\sum_{k=0}^{n}{n \brace k}_{\lambda}X^{k}D^{k}e^{x}\label{20}\\
&=\frac{1}{e^{x}}\sum_{k=0}^{n}{n \brace k}_{\lambda}x^{k}e^{x}\nonumber\\
&=\phi_{n,\lambda}(x),\quad (n\ge 0). \nonumber
\end{align}
Therefore, by \eqref{20}, we obtain the following theorem.
\begin{theorem}
For $n\ge 0$, we have
\begin{equation}
\frac{1}{e^{x}}\big(XD\big)_{n,\lambda}e^{x}=\phi_{n,\lambda}(x). 	\label{21}
\end{equation}
\end{theorem}
By \eqref{2}, we get
\begin{align}
\big(XD\big)_{m+n,\lambda}&=\big(XD\big)\cdots\big(XD-(m-1)\lambda\big)\big(XD-m\lambda\big)\cdots\big(XD-m\lambda-(n-1)\lambda\big)\label{22}\\
&=\big(XD\big)_{m,\lambda}\big(XD-m\lambda\big)_{n,\lambda}=\big(XD-m\lambda\big)_{n,\lambda}\big(XD\big)_{m,\lambda},\nonumber
\end{align}
where $m,n$ are nonnegative integers. \par
We note that
\begin{equation}
(x+y)_{n,\lambda}=\sum_{k=0}^{n}\binom{n}{k}(x)_{k,\lambda}(y)_{n-k,\lambda},\quad (n\ge 0). \label{23}
\end{equation} \\
For $n\ge 0$, by \eqref{17} and \eqref{23}, we get
\begin{align}
&\big(XD-m\lambda\big)_{n,\lambda}X^{j}=\big(XD-m\lambda\big)\big(XD-m\lambda-\lambda\big)\cdots\big(XD-m\lambda-(n-1)\lambda\big)X^{j}	\label{24} \\
&=\big(XD-m\lambda\big)\cdots\big(XD-m\lambda-(n-2)\lambda\big)X^{j}\big(XD-m\lambda-(n-1)\lambda+j\big)\nonumber\\
&=\cdots\nonumber\\
&=X^{j}\big(XD+j-m\lambda\big)\big(XD+j-m\lambda-\lambda\big)\cdots\big(XD+j-m\lambda-(n-1)\lambda\big)\nonumber\\
&=X^{j}\big(XD+j-m\lambda\big)_{n,\lambda}=X^{j}\sum_{k=0}^{n}\binom{n}{k}(XD)_{k,\lambda}(j-m\lambda)_{n-k,\lambda}.\nonumber
\end{align} \par
From \eqref{19}, \eqref{22} and \eqref{24}, we note that
\begin{align}
\big(XD)_{m+n,\lambda}&=\big(XD-m\lambda\big)_{n,\lambda}\big(XD\big)_{m,\lambda}=\sum_{j=0}^{m}{m \brace j}_{\lambda}\big(XD-m\lambda\big)_{n,\lambda}X^{j}D^{j} \label{25}\\
&=\sum_{j=0}^{m}\sum_{k=0}^{n}{m \brace j}_{\lambda}\binom{n}{k}\big(j-m\lambda\big)_{n-k,\lambda}X^{j}\big(XD\big)_{k,\lambda}D^{j}, \nonumber	
\end{align}
where $m,n$ are nonnegative integers. \\
For $m,n\ge 0$, by \eqref{21} and \eqref{25}, we get
\begin{align}
\phi_{m+n,\lambda}(x)&=\frac{1}{e^{x}}\big(XD\big)_{m+n,\lambda}e^{x} 	\label{26}\\
&=\sum_{j=0}^{m}\sum_{k=0}^{n}{m \brace j}_{\lambda}\binom{n}{k}(j-m\lambda)_{n-k,\lambda}\frac{1}{e^{x}}X^{j}\big(XD\big)_{k,\lambda}D^{j}e^{x}\nonumber\\
&=\sum_{j=0}^{m}\sum_{k=0}^{n}{m \brace j}_{\lambda}\binom{n}{k}(j-m\lambda)_{n-k,\lambda}\frac{1}{e^{x}}X^{j}\big(XD\big)_{k,\lambda}e^{x}\nonumber\\
&=\sum_{j=0}^{m}\sum_{k=0}^{n}{m \brace j}_{\lambda}\binom{n}{k}(j-m\lambda)_{n-k,\lambda}x^{j}\phi_{k,\lambda}(x).\nonumber
\end{align}
Therefore, by \eqref{26}, we obtain the following theorem.
\begin{theorem}
For $m,n\ge 0$, we have
\begin{displaymath}
\phi_{m+n,\lambda}(x)= \sum_{j=0}^{m}\sum_{k=0}^{n}\binom{n}{k}{m \brace j}_{\lambda}(j-m\lambda)_{n-k,\lambda}x^{j}\phi_{k,\lambda}(x).
\end{displaymath}
Letting $x=1$ yields
\begin{equation}
\phi_{m+n,\lambda}=\sum_{j=0}^{m}\sum_{k=0}^{n}\binom{n}{k}{m \brace j}_{\lambda}(j-m\lambda)_{n-k,\lambda}\phi_{k,\lambda}. \label{26-1}
\end{equation}
\end{theorem}
We observe that taking $\lambda \rightarrow 0$ in \eqref{26-1} gives the Spivey's formula for Bell numbers in \eqref{0}. \par
Let $f(x)=x^{m},\ (m \ge 0)$. For $n,r\ge 0$, by \eqref{6}, we get
\begin{align}
\big(XD+r\big)_{n,\lambda}f(x)&=\big(XD+r\big)\big(XD+r-\lambda\big)\cdots\big(XD+r-(n-1)\lambda\big)x^{m}\label{27}\\
&=(m+r)_{n,\lambda}x^{m}=\sum_{k=0}^{n}{n+r \brace k+r}_{r,\lambda}(m)_{k}x^{m}\nonumber\\
&=\sum_{k=0}^{n}{n+r \brace k+r}_{r,\lambda}X^{k}D^{k}x^{m}=\sum_{k=0}^{n}{n+r \brace k+r}_{r,\lambda}X^{k}D^{k}f(x).\nonumber	
\end{align}
Thus, by \eqref{27}, we get
\begin{equation}
\big(XD+r\big)_{n,\lambda}=\sum_{k=0}^{n}{n+r \brace k+r}_{r,\lambda}X^{k}D^{k},\quad (n,r \ge 0).\label{29}
\end{equation} \\
For $n,r\ge 0$, by \eqref{12} and \eqref{14}, we have
\begin{align}
\frac{1}{e^{x}}\big(XD+r\big)_{n,\lambda}e^{x}&=\frac{1}{e^{x}}\sum_{k=0}^{n}{n+r \brace k+r}_{r,\lambda}X^{k}D^{k}e^{x}\label{30}\\
&=\sum_{k=0}^{n}{n+r \brace k+r}_{r,\lambda}x^{k}=\phi_{n,\lambda}^{(r)}(x). \nonumber
\end{align}
Therefore, by \eqref{30}, we obtain the following theorem.
\begin{theorem}
For $n,r\ge 0$, we have
\begin{equation}
\frac{1}{e^{x}}\big(XD+r\big)_{n,\lambda}e^{x}	=\phi_{n,\lambda}^{(r)}(x).\label{31}
\end{equation}
\end{theorem}
For $m,n\ge 0$, by \eqref{23} and \eqref{29}, we get
\begin{align}
\big(XD+r\big)_{n+m,\lambda}&=\big(XD+r-m\lambda\big)_{n,\lambda}\big(XD+r\big)_{m,\lambda}  \label{32}\\
&=\big(XD+r-m\lambda\big)_{n,\lambda}\sum_{k=0}^{m}{m+r \brace k+r}_{r,\lambda}X^{k}D^{k}\nonumber\\
&=\sum_{k=0}^{m}{m+r \brace k+r}_{r,\lambda}\big(XD+r-m\lambda\big)_{n,\lambda}X^{k}D^{k}.\nonumber
\end{align}
Now, we observe that
\begin{align}
\big(XD+r-m\lambda\big)_{n,\lambda}X^{k}
&=X^{k}\big(XD+r+k-m\lambda\big)_{n,\lambda} \label{33}\\
&=X^{k}\sum_{l=0}^{n}\binom{n}{l}\big(XD+r\big)_{l,\lambda}(k-m\lambda)_{n-l,\lambda}.  \nonumber
\end{align}
By \eqref{32} and \eqref{33}, we get
\begin{equation}
\big(XD+r\big)_{n+m,\lambda}=\sum_{k=0}^{m}\sum_{l=0}^{n}{m+r \brace k+r}_{r,\lambda}\binom{n}{l}(k-m\lambda)_{n-l,\lambda}X^{k}\big(XD+r\big)_{l,\lambda}D^{k}.\label{34}
\end{equation}
From \eqref{31} and \eqref{34}, we note that
\begin{align}
\phi_{n+m,\lambda}^{(r)}(x)&=\frac{1}{e^{x}}\big(XD+r\big)_{n+m,\lambda}e^{x}\label{35}\\
&=\sum_{k=0}^{m}\sum_{l=0}^{n}{m+r \brace k+r}_{r,\lambda}\binom{n}{l}\big(k-m\lambda\big)_{n-l,\lambda}\frac{1}{e^{x}}X^{k}\big(XD+r\big)_{l,\lambda}D^{k}e^{x}\nonumber\\
&=\sum_{k=0}^{m}\sum_{l=0}^{n}{m+r \brace k+r}_{r,\lambda}\binom{n}{l}(k-m\lambda)_{n-l,\lambda}x^{k}\phi_{l,\lambda}^{(r)}(x),\nonumber	
\end{align}
where $n,m$ are nonnegative integers. \\
Therefore, by \eqref{35}, we obtain the following theorem.
\begin{theorem}
For $n,m\in\mathbb{Z}$ with $n,m\ge 0$, we have
\begin{equation*}
\phi_{n+m,\lambda}^{(r)}(x)=\sum_{k=0}^{m}\sum_{l=0}^{n}\binom{n}{l}{m+r \brace k+r}_{r,\lambda}(k-m\lambda)_{n-l,\lambda}x^{k}\phi_{l,\lambda}^{(r)}(x).
\end{equation*}
Letting $x=1$ gives
\begin{equation*}
\phi_{n+m,\lambda}^{(r)}=\sum_{k=0}^{m}\sum_{l=0}^{n}\binom{n}{l}{m+r \brace k+r}_{r,\lambda}(k-m\lambda)_{n-l,\lambda}\phi_{l,\lambda}^{(r)}.
\end{equation*}
\end{theorem}
Note that
\begin{displaymath}
\phi_{n+m}^{(r)}=\lim_{\lambda\rightarrow 0}\phi_{n+m,\lambda}^{(r)}=\sum_{k=0}^{m}\sum_{l=0}^{n}\binom{n}{l}{m+r \brace k+r}_{r}k^{n-l}\phi_{l}^{(r)}.
\end{displaymath}

\section{Conclusion}
In this paper, we derived Spivey's type recurrence relations for the degenerate Bell polynomials $\phi_{n,\lambda}(x)$ and the degenerate $r$-Bell polynomials $\phi_{n,\lambda}^{(r)}(x)$ by using the `multiplication by $x$' operator $X$ and the differentiation operator $D=\frac{d}{dx}$. \par
Here we review the idea of our proof for $\phi_{n,\lambda}(x)$. By using \eqref{22}, \eqref{19}, \eqref{17} and \eqref{23} in this order, we can proceed as follows:
\begin{align}
(XD)_{m+n,\lambda}&=(XD-m \lambda)_{n,\lambda}(XD)_{m,\lambda} \label{36}\\
&=(XD-m \lambda)_{n,\lambda}\sum_{j=0}^{m}{m \brace j}_{\lambda} X^{j}D^{j}\nonumber \\
&=\sum_{j=0}^{m}{m \brace j}_{\lambda}(XD-m \lambda)_{n,\lambda}X^{j}D^{j}\nonumber \\
&=\sum_{j=0}^{m}{m \brace j}_{\lambda}X^{j}(XD+j-m \lambda)_{n,\lambda}D^{j} \nonumber\\
&=\sum_{j=0}^{m}{m \brace j}_{\lambda}X^{j}\sum_{k=0}^{n}\binom{n}{k}(j-m \lambda)_{n-k, \lambda}(XD)_{k,\lambda} D^{j} \nonumber \\
&=\sum_{j=0}^{m}\sum_{k=0}^{n}\binom{n}{k}{m \brace j}_{\lambda}(j-m \lambda)_{n-k, \lambda}X^{j}(XD)_{k,\lambda} D^{j}. \nonumber
\end{align}
Using \eqref{36} and $\phi_{n,\lambda}(x)=\frac{1}{e^{x}}(XD)_{n,\lambda}e^{x}$ (see \eqref{21}), we obtain the recurrence relation in Theorem 2.2:
\begin{equation*}
\phi_{m+n,\lambda}(x)= \sum_{j=0}^{m}\sum_{k=0}^{n}\binom{n}{k}{m \brace j}_{\lambda}(j-m\lambda)_{n-k,\lambda}x^{j}\phi_{k,\lambda}(x).
\end{equation*}

\end{document}